\documentclass[12pt]{amsart}

\usepackage{epsfig}
\usepackage{verbatim}

\begin{document}
\title {Matroids with at least two regular elements}

\maketitle
\begin {center}
S. R. Kingan 
\footnote{The first author is partially supported by a PSC-CUNY Award 63076-00 41.} \\                                                                                                    
Department of Mathematics \\
Brooklyn College, City University of New York\\
 Brooklyn, NY 11210\\
skingan@brooklyn.cuny.edu\\  
\end {center}

\begin {center}
Manoel Lemos 
\footnote{The second author is partially supported by CNPq under grant number 300242/2008-05.}\\
Departamento de Matematica \\
Universidade Federal de Pernambuco\\
Recife, Pernambuco, 50740-540, Brazil\\
manoel@dmat.ufpe.br\\  
\end {center}
\bigskip

\begin{abstract} For a matroid $M$, an element $e$ such that both $M\backslash e$ and $M/e$ are regular is called a regular element of $M$. We determine completely the structure of non-regular matroids with at least two regular elements. Besides four small size matroids, all 3-connected matroids in the class can be pieced together from $F_7$ or $S_8$ and a regular matroid using 3-sums. This result takes a step toward solving a problem posed by Paul Seymour: find all 3-connected non-regular matroids with at least one regular element [5, 14.8.8].
\end{abstract}

\bigskip 

\section {Introduction}

The matroid terminology follows Oxley [5]. Let $M$ be a matroid and $X$ be a subset of the ground set $E$. The {\it connectivity function}  $\lambda$ is defined as $\lambda (X) = r(X) + r(E-X) - r(M)$.  Observe that $\lambda (X) = \lambda (E-X)$.  For $j\ge 1$, a partition $(X_1, X_2)$ of $E$ is called a $j$-separation if $|X_1|, |X_2|\ge j$, and $\lambda (X_1) \le j-1$.  When $\lambda (X_1)=j-1$, we call $(X_1, X_2)$ an {\it exact j-separation}.  When $\lambda (X_1)=j-1$ and $|X_1|=j$ or $|X_2|=j$ we call $(X_1, X_2)$ a {\it minimal exact j-separation}. For $k\ge 2$, we say $M$ is {\it k-connected} if $M$ has no $j$-separation for $j\le k-1$.  A matroid is {\it internally $k$-connected} if it is $k$-connected and has no non-minimal exact $k$-separations. In particular, a simple matroid is 3-connected if $\lambda (X_1)\ge 2$ for all partitions $(X_1, X_2)$ with $|X_1|, |X_2|\ge 3$. A 3-connected matroid is {\it internally $4$-connected}  if $\lambda (X_1)\ge 3$ for all partitions $(X_1, X_2)$ with $|X_1|, |X_2|\ge 4$. 

The 1-sum, 2-sum, and 3-sum of binary matroids are defined in [6]. A {\it cycle} of a binary matroid is a disjoint union of circuits. Let $M_1$ and $M_2$ be binary matroids with non-empty ground sets $E_1$ and $E_2$, respectively. We define a new binary matroid $M_1 \triangle M_2$ to be the matroid with ground set $E_1\triangle E_2$ and with cycles having the form $C_1\triangle C_2$ where $C_i$ is a cycle of $M_i$ for $i=1, 2$. When $E_1\cup E_2=\phi$, then $M_1 \triangle M_2$ is a {\it $1$-sum} of $M_1$ and $M_2$. When $|E_1|, |E_2|\ge 3$, $E_1\cap E_2=\{z\}$ and $z$ is not a loop or coloop of $M_1$ or $M_2$, then $M_1 \triangle M_2$ is a {\it $2$-sum} of $M_1$ and $M_2$. When $|E_1|, |E_2|\ge 7$, $E_1\cap E_2=T$ and $T$ is a triangle in $M_1$ and $M_2$, then $M_1 \triangle M_2$ is a {\it $3$-sum} of $M_1$ and $M_2$.

An element $e$ in a non-regular matroid $M$ is called a {\it regular element} if both $M\backslash e$ and $M/e$ are regular. Seymour posed the following problem  that appears in Oxley's book {\it Matroid Theory} [5, 14.8.8]:  Find all 3-connected non-regular matroids with at least one regular element. In other words the problem is to find all 3-connected non-regular elements with at least one regular element.  In this paper, we take a step toward solving this problem by determining the class of non-regular matroids with at least two regular elements.

We denote the 4-point line as $U_{2, 4}$ and the  Fano matroid as $F_7$. We denote by $S_8$ the following single-element extension of $F_7$. It is self-dual. A single-element extension of $S_8$ that will play a role is $P_9$ shown below. 

\small
\[ 
F_7=\left[ 
\begin{array}{ccc|cccc}
&&&   1&0&1&1 \\
&I_3&&1&1&0&1 \\
&&&   0&1&1&1 \\
\end{array} 
\right]  
S_8=\left[ 
\begin{array}{ccc|cccc}
&&&   0&1&1&1 \\
&I_4&&1&0&1&1 \\
&&&   1&1&0&1 \\
&&&   1&1&1&1
\end{array} 
\right]  
P_9= \left[ 
\begin{array}{ccc|ccccc}
&&&   0&1&1&1&1 \\
&I_4&&1&0&1&1&1 \\
&&&   1&1&0&1&0 \\
&&&   1&1&1&1&0
\end{array} 
\right]  
\]
\normalsize

\noindent 
 For this paper, it helps to think of $F_7$ as the single-element extension of the 3-wheel with spokes labeled $\{1, 2, 3\}$ where the new element forms a circuit with $\{1, 2, 3\}$. The matroid $P_9$ is the single-element extension of the 4-wheel with spokes $\{1, 2, 3, 4\}$ where the new element forms a circuit with any three consecutive spokes, say $\{1, 2, 3\}$. Then  $P_9\backslash 1\cong S_8$ and $P_9\backslash 3\cong S_8$. Moreover, $P_9\backslash \{1, 3\}\cong F_7^*$.

 Let $F_7^p$ and $S_8^p$ be the matroids obtained from $F_7$ and $S_8$, respectively, by adding an element in parallel with an element belonging to at least two triangles. Note that every element of $F_7$ is in at least two triangles, but only one element of $S_8$ is in two triangles. The main result of this paper gives a complete characterization of the matroids with at least two regular elements.

\bigskip

\noindent {\bf Theorem 1.1.} {\it A $3$-connected non-regular matroid $M$ has at least two regular elements if and only if
\begin{enumerate}
\item[(i)] $M$ is  $U_{2,4}$, $F_7,F_7^*$ or $S_8$; or
\item[(ii)] $M$ is the $3$-sum of $F_7$ or $S_8$ with a $3$-connected regular matroid (with the possible exception of elements in parallel with the $3$-sum triangle); or
\item[(iii)] $M$ is the $3$-sum of $F_7^p$ or $S_8^p$ with two $3$-connected regular matroids (with the possible exception of elements in parallel with the $3$-sum triangle). These two $3$-sums are made along two disjoint triangles of $F_7^p$ or $S_8^p$.
\end{enumerate}}
\bigskip

In order to prove this result we use the following theorems. The first is by Oxley and appears in [4, 3.9]:

\bigskip

\noindent {\bf Theorem 1.2.} {\it Let $M$ be a $3$-connected matroid having an element $e$ such that $M\backslash e$ and $M/e$ are both regular. Then $M\cong U_{2,4}$.} $\qed$
\bigskip

\noindent The next result by Zhou appears in [7, 1.2]. The matroid $S_{10}$, shown below, is the first matroid in the internally 4-connected infinite family of almost-graphic matroids $S_{3n+1}$ [3].  The matroid $M(E_5)$ appears in [1] where Kingan characterized the class of matroids with no minors isomorphic to $M(K_5\backslash e)$, $M^*(K_5\backslash e)$ and $AG(3, 2)$. $M(E_5)$ is a splitter for this class.  It  is self dual and internally 4-connected. The self-dual 4-connected matroid $T_{12}$ appears in [2].  

\tiny
\[ 
S_{10}=\left[ 
\begin{array}{ccc|cccccc}
&&&   1&0&0&1&1&0 \\
&I_4&&1&1&0&0&1&1 \\
&&&   0&1&1&0&1&1 \\
&&&   0&0&1&1&0&1
\end{array} 
\right]  
E_5=\left[ 
\begin{array}{ccc|cccccc}
&&&    0&1&1&1&1 \\
&&&    1&0&1&1&0 \\
&I_5&& 1&1&0&1&1 \\
&&&    1&1&1&1&0 \\
&&&    1&1&0&0&0
\end{array} 
\right] 
T_{12}=\left[ 
\begin{array}{ccc|cccccc}
&&&   1&1&0&0&0&1 \\
&&&   1&0&0&0&1&1 \\
&I_6&&0&0&0&1&1&1 \\
&&&   0&0&1&1&1&0 \\
&&&   0&1&1&1&0&0 \\
&&&  1&1&1&0&0&0 
\end{array} 
\right] 
\] 
\normalsize

\bigskip

\noindent {\bf Theorem 1.3.}  {\it A non-regular internally $4$-connected binary matroid other than $F_7$ and $F_7^*$ contains one of the following matroids as a minor: $M(E_5)$, $S_{10}$, $S_{10}^*$, $T_{12}\backslash e$, and $T_{12}/e$.} $\qed$
\bigskip

\noindent Finally, we use the following results by Seymour that appear in [6, 2.9] and [6, 4.1]:
\bigskip

\noindent {\bf Theorem 1.4.}  {\it  If $(X_1, X_2)$ is an exact $3$-separation of a binary matroid $M$, with $|X_1|, |X_2|\ge 4$, then there are binary matroids $M_1$, $M_2$ on $X_1\cup T$, $X_2\cup T$, respectively (where $T$ contains three new elements), such that $M$ is the $3$-sum of $M_1$ and $M_2$. Conversly if $M$ is the $3$-sum of $M_1$ and $M_2$, then $(E(M_1)-E(M_2), E(M_2)-E(M_1))$ is an exact $3$-separation of $M$, and $|E(M_1)-E(M_2)|, |E(M_2)-E(M_1)|\ge 4$.} $\qed$
\bigskip

\noindent {\bf Theorem 1.5.}  {\it  If $M$ is binary and is the $3$-sum of $M_1$ and $M_2$, and $M$ is $3$-connected, then $M_1$ and $M_2$ are isomorphic to minors of $M$.} $\qed$
\bigskip

\noindent In the next section we give several separation lemmas that are used in the proof of the main theorem. In the third section we give results on the number of regular elements in a matroid. Finally, in the fourth section we prove Theorem 1.1.

\section{Understanding 3-separations in the context of regular elements}

Let $M$ be a 3-connected non-regular binary matroid such that $M$ is the 3-sum of matroids $M_1$ and $M_2$ where $|E(M_1)|, |E(M_2)|\ge 7$, $E(M_1)\cap E(M_2)=T$ and $T$ is a triangle in $M_1$ and $M_2$.  Assume that $e\in E(M_1)-E(M_2)$ is a regular element of $M$.  
\bigskip

\noindent {\bf Lemma 2.1.} {\it The element $e$ is not spanned by $E(M_2)-E(M_1)$ in $M$.}
\bigskip

\noindent {\bf Proof.} Suppose $e$ is spanned by $E(M_2)-E(M_1)$ in $M$. Then $e$ is spanned by $T$ in $M_1$ and so $e$ is in parallel to some element $t\in T$. By hypothesis, $M\backslash e$ is regular. Observe that $M_1\backslash e$ and $M_2$ are regular because:
\begin{enumerate}
\item[(i)] when $|E(M_1)|>7$, $M\backslash e$ is the 3-sum of $M_1\backslash e$ with $M_2$; and
\item[(ii)] when $|E(M_1)|=7$, $M_1\backslash e$ has 6 elements and is isomorphic to $M(K_4)$. So $M\backslash e$ is obtained from $M_2$ after a $\Delta-Y$ operation along the triangle $T$.
\end{enumerate}
But $M_1$ is obtained from $M_1\backslash e$ by adding $e$ in parallel with $t$. Therefore $M_1$ and $M_2$ are regular; a contradiction because the class of regular matroids is closed under 3-sums. Thus $e$ is not spanned by $E(M_2)-E(M_1)$ in $M$. $\qed$
\bigskip

\noindent {\bf Lemma 2.2.} {\it The element $e$ is not spanned by $E(M_2)-E(M_1)$ in $M^*$.}
\bigskip

\noindent {\bf Proof.} If $N_i$ is obtained from $M_i$ by a $\Delta - Y$ operation along the triangle $T$, then $M^*$ is the 3-sum of $N_1^*$ and $N_2^*$. Applying Lemma 2.1, we conclude that $e$ is not spanned by $E(M_2)-E(M_1)$ in $M^*$.  $\qed$
\bigskip

 In the next result, we describe how the presence of a regular element in $M_1$ impacts the structure of $M$. We prove that one of two situations occur: either $M_1$ is non-regular with $e$ as a regular element and $M_2$ is regular or  $M_2$ is non-regular and $M_1$ is a small matroid with with a specific structure.  In the latter situation we prove that $E(M_1)-T=T'\cup T^*$ where $T'$ is a triangle and $T^*$ is a triad such that $e\in T\cap T^*$ and $E(M_1)-E(M_2)$ is closed in $M$. Since $M$ is binary, a triangle and triad must intersect in an even number of elements. This means $M_1$ has just 7 elements, one of which is parallel with an element of $T$.

\bigskip

\noindent {\bf Lemma 2.3.} {\it 
\begin{enumerate}
\item[(i)] $M_2$ is a regular matroid; or
\item[(ii)] there is a triangle $T'$ and a triad $T^*$ of $M$ such that $e\in T'\cap T^*$ and $E(M_1)-T=T'\cup T^*$.
\end{enumerate}
Moreover,
\begin{enumerate}
\item[(iii)] when (i) happens, $M_1$ is a non-regular matroid having $e$ as a regular element; and
\item[(iv)] when (ii) happens, $E(M_1)-E(M_2)$ is closed in $M$.
\end{enumerate}}
\bigskip

\noindent {\bf Proof.} Assume that (i) does not hold, that is,
\begin{equation}\label{27.01.11.a}
\hbox{$M_2$ is non-regular.}
\end{equation}
First, we establish that
\begin{equation}\label{26.01.11.a}
r(M_1)=3\hbox{ or ${\rm si}(M/e)$ is not 3-connected.}
\end{equation}
Suppose that $r(M_1)\ge 4$ and ${\rm si}(M/e)$ is 3-connected. If $T'$ is a triangle of $M$ containing $e$, then, by Lemma 2.1, $|E(M_2)\cap T'|\le 1$. Therefore we may assume that ${\rm si}(M/e)=M/e\backslash X$, for $X\subseteq E(M_1)-T$. If $M_1/e\backslash X\simeq M(K_4)$, then $M_2$ is obtained from ${\rm si}(M/e)$ after a $Y-\Delta$ operation along the triad $E(M_1)-(e\cup X\cup T)$. So $M_2$ is regular; a contradiction to~\eqref{27.01.11.a}. If $M_1/e\backslash X\not\simeq M(K_4)$, then ${\rm si}(M/e)$ is the 3-sum of $M_1/e\backslash X$ and $M_2$. As ${\rm si}(M/e)$ is regular, it follows that $M_2$ is regular; a contradiction to~\eqref{27.01.11.a}. We have~\eqref{26.01.11.a}.

If $N_i$ is obtained from $M_i$ by a $\Delta - Y$ operation along the triangle $T$, then $M^*$ is the 3-sum of $N_1^*$ and $N_2^*$. Note that Lemma 2.3(i) holds for the decomposition $M=M_1\bigtriangleup M_2$ if and only if Lemma 2.3(i) holds for the decomposition $M^*=N_1^*\bigtriangleup N_2^*$. The analogous statement occurs when we replace (i) by (ii). Therefore, the dual of~\eqref{26.01.11.a} becames
\begin{equation}\label{07.02.11.a}
r(N_1^*)=3\hbox{ or $[{\rm co}(M\backslash e)]^*={\rm si}(M^*/e)$ is not 3-connected.}
\end{equation}
By Bixby's Theorem [5, 8.4.6], ${\rm si}(M/e)$ or ${\rm co}(M\backslash e)$ is 3-connected. By~\eqref{26.01.11.a} and~\eqref{07.02.11.a}, $r(M_1)=3$ or $r(N_1^*)=3$. Taking the dual when necessary, we may assume that
\begin{equation}\label{07.02.11.b}
r(M_1)=3.
\end{equation}
Next, we prove the following claim:
\bigskip

\noindent {\bf Claim:} $M_1$ does not have a minor $N$ such that $T$ and $T'=E(N)-T$ are triangles of $N$, $e\not\in E(N)=T\cup T'$ and $r(N)=2$.
\bigskip

Suppose that $N$ exists, say $N=M_1\backslash X/Y$. By hypothesis, $e\in X\cup Y$ and so $M\backslash X/Y$ is regular. Moreover, $M\backslash X/Y$ is isomorphic to $M_2$. Thus $M_2$ is regular; a contradiction to~\eqref{27.01.11.a}. Therefore the claim holds.
\bigskip

If ${\rm si}(M_1)\simeq F_7$, then $M_1/e$ is a rank-2 matroid. By Lemma 2.1, $M_1/e$ has $T$ as a triangle. We have a contradiction by the Claim because every parallel class of $M_1/e$ is non-trivial. Hence, by~\eqref{07.02.11.b}, ${\rm si}(M_1)\simeq M(K_4)$. In particular, $T^*=E(M_1)-{\rm cl}_{M_1}(T)$ is a triad of $M_1$. By Lemma 2.1, $e\in T^*$, say $T^*=\{e,e_1,e_2\}$. Let $f_1,\dots,f_k$ be the elements of ${\rm cl}_{M_1}(T)-T$. For each $i$, there is $t_i\in T$ such that $\{f_i,t_i\}$ is a parallel class of $M_1$. By the Claim, $k\le 2$. Next, we establish that
\begin{equation}\label{26.01.11.b}
k=1.
\end{equation}
As $|E(M_1)|\ge 7$ and $|E(M_1)-{\rm cl}_{M_1}(T)|=3$, it follows that $k\ge 1$. If~\eqref{26.01.11.b} does not hold, then $k=2$. In $M_1/e$, by the Claim, $e_i$ is in parallel with $f_j$, say $e_i$ is in parallel with $f_i$, for both $i$. Therefore $T_i=\{e,e_i,f_i\}$ is a triangle of $M$, for both $i$, and so $T_1\bigtriangleup T_2\bigtriangleup\{f_1,f_2,t_3\}=\{e_1,e_2,t_3\}$, where $T=\{t_1,t_2,t_3\}$ is a triangle of $M_1$. Thus $N=M_1\backslash e/e_1$ is a minor of $M_1$ contrary to the Claim. Thus~\eqref{26.01.11.b} holds. By the claim $e_1$ or $e_2$ is in parallel with $f_1$ in $M_1/e$, say $e_1$. That is, $T'=\{e,e_1,f_1\}$ is a triangle of $M_1$ and so of $M$. We have (ii).
 
Assume that (i) happens, that is, $M_2$ is regular. Thus $M_1$ is non-regular because $M$ is non-regular. To conclude (iii) we need to prove only that $e$ is a regular element of $M_1$. By the proof of Theorem 1.5, there are disjoint subsets $Y$ and $Z$ of $E(M_2)-E(M_1)$ such that $N=M_2\backslash Y/Z$ is a 6-element matroid such that $T''=E(N)-T$ is a triangle of $N$ and, for each $f\in T$, there is an $f''\in T''$ such that $\{f,f''\}$ is a circuit of $N$. So $M\backslash Y/Z$ is isomorphic to $M_1$ --- this isomorphism fix each element of $E(M_1)-E(M_2)$ and sends $f''$ into $f$, for each $f''\in T''$. As both $M\backslash e$ and $M/e$ are regular, it follows that $(M\backslash e)\backslash Y/Z\simeq M_1\backslash e$ and $(M/e)\backslash Y/Z\simeq M_1/e$ are regular. That is, $e$ is a regular element of $M_1$. We have (iii).

Assume that (ii) happens. If $E(M_1)-E(M_2)$ spans an element $g$ of $E(M_2)-E(M_1)$ in $M$, then $[E(M_1)-E(M_2)]\cup g$ is a 3-separating set for $M$. Using the 3-separation induced by this set, we can decompose $M$ as the 3-sum of matroids $M_1'$ and $M_2'$ such that $E(M_1')=[E(M_1)-E(M_2)]\cup g\cup T''$ and $T''=E(M_1')\cap E(M_2')$. Note that, in $M_1'$, the element $g$ is in parallel with some element of $T''$. In particular, $M_1'\backslash g\simeq M_1$ is regular. So $M_1'$ is regular; a contradiction to this lemma. Thus $E(M_1)-E(M_2)$ is closed in $M$. $\qed$
\bigskip

Now that we have shown  $M$ has a clearly defined structure, we want to say more about the second situation.  Recall that $R(M)$ is the set of regular elements. 
For a triangle $T'$ and triad $T^*$ of $M$, we say that $T',T^*$ is an {\it undesired fan} if $T'\cap T^*\cap R(M)\not=\emptyset$. Note that $\{T'\cup T^*, E(M)-(T'\cup T^*)\}$ is an exact 3-separation for $M$ and by Theorem 1.4, it is possible to decompose $M$ as a 3-sum using it. In the next lemma we show that the presence of an undesired fan implies the existence of two regular elements.

\bigskip

\noindent {\bf Lemma 2.4.} {\it If $T',T^*$ is an undesired fan in $M$ such that $E(M_1)-E(M_2)=T'\cup T^*$, then $T'\cap T^*\subseteq R(M)$. Moreover, if $T^*-T'=\{f\}$, then $M/f$ is a $3$-connected non-regular matroid such that $T'\cap T^*\subseteq R(M/f)$.}
\bigskip

\noindent {\bf Proof.} Suppose that $T'=\{e,e',t\}$, $T^*=\{e,e',f\}$ and $e\in R(M)$. In $M/e'$, $t$ and $e$ are in parallel. As $M\backslash e$ and so $M/e'\backslash e$ is regular, it follows that $M/e'$ is regular because $M/e'$ is obtained from $M/e'\backslash e$ by adding $e$ in parallel with $t$. Using duality, we conclude that $M\backslash e'$ is regular. Hence $e'$ is a regular element of $M$ and so $T'\cap T^*\subseteq R(M)$. 

Next,  observe that $E(M_1)=T'\cup T^*\cup T$ and $E(M_2)=[E(M)-(T'\cup T^*)]\cup T$. As $M_1$ is regular, it follows that $M_2$ is non-regular. By Lemma 2.3, $f$ does not belong to a triangle of $M$. So $M/f$ is 3-connected because ${\rm si}(M/f)$ is 3-connected. But $M/f\simeq M_2$ because $M_1/f$ has three non-trivial parallel classes each containing one element of $T'$ and another of $T$. The result follows because $R(M)\subseteq R(M/f)$. $\qed$
\bigskip

In the next lemma, we prove that, when this happens, it is possible to uncontract $f$ keeping the property of these two regular elements.
\bigskip

\noindent {\bf Lemma 2.5.} {\it Let $N$ be a $3$-connected non-regular binary matroid having different regular elements $e$ and $e'$. Suppose that $T'$ is a triangle of $N$ such that $e,e'\in T'$ and $\{e,e'\}$ is not contained in a triad of $N$. If $M$ is a one-element binary lift of $N$, say $M/f=N$, such that $\{e,e',f\}$ is a triad of $M$, then $e$ and $e'$ are regular elements of $M$ (and $M$ is $3$-connected).}
\bigskip

\noindent {\bf Proof.} Observe that ${\rm si}(M/e)=M/e\backslash e'$. But, in $M\backslash e'$, $e$ and $f$ are in series. So $M/e\backslash e'\simeq M/f\backslash e'=N\backslash e'$ and  ${\rm si}(M/e)$ is regular. Thus $M/e$ is regular. As $M\backslash e/f=N\backslash e$, it follows that $M\backslash e/f$ is regular and so $M\backslash e$ is regular. That is, $e$ is a regular element of $M$. A similar argument holds with $e'$. $\qed$
\bigskip

\section {The number of regular elements in a matroid}

Next we prove a result on the number of regular elements in a binary non-regular matroid. Observe that, $F_7^*$ has two single-element extensions $S_8$ and $AG(3, 2)$. The matroid $AG(3, 2)$ has one single-element extension $Z_4$. The matroid $S_8$ has two single-element extensions, $Z_4$ and $P_9$. Observe further that $F_7$ and $F_7^*$ have seven regular elements and $P_9$ has four regular elements. $AG(3, 2)$ has zero regular elements and consequently so do all its 3-connected extensions and coextensions.

\bigskip

\noindent {\bf Lemma 3.1.} {\it Let $M$ be a $3$-connected non-regular binary matroid. If $|E(M)|\ge 9$, then $|R(M)|=0, 1, 2$ or $4$. Moreover, if $|R(M)|=4$, then $R(M)$ is both a circuit and a cocircuit of $M$.}
\bigskip

\noindent {\bf Proof.} Since $|E(M)|\ge 9$, $M$ must have $P_9$ or $P_9^*$ as a minor and since $|R(P_9)|=4$, it follows that $|R(M)|\le 4$. 

Suppose $|R(M)|=3$. Choose a minimal counter-example $M$. First observe that $|E(M)|\ge 10$ and by Theorem 1.3 $M$ is not internally 4-connected. By Theorem 1.4, we can decompose $M$ as the 3-sum of matroids $M_1$ and $M_2$ such that $E(M_1)\cap E(M_2)=T$ and $E(M_1)\cap R(M)\neq \phi$. If Lemma 2.3 (ii) occurs and $f\in T^*-T'$, then by Lemma 2.4 and the choice of $M$, $M/f$ has $R(M)\cup g$ as a circuit-cocircuit, for some element $g$. Note that $[R(M)\cup g]\bigtriangleup T^*$ is a triad of $M$ and $[R(M)\cup g]\bigtriangleup T'$ is a triangle of $M$ whose intersection contains a regular element. Therefore, by Lemma 2.4 the intersection has two regular elements ($g$ is the other regular element); a contradiction. Thus Lemma 2.3(i) occurs. Observe that $R(M)$ is contained in a circuit-cocircuit of $M_1$ consisting of regular elements avoiding $T$. Thus every element in this circuit-cocircuit is also a regular element of $M$; a contradiction. Thus we proved that $M$ cannot have exactly three regular elements. 

Next, suppose $|R(M)|=4$, but $R(M)$ is not a circuit cocircuit. Choose a minimum counterexample. As before, $|E(M)|\ge 10$ and by Theorem 1.3 $M$ is not internally 4-connected. By Theorem 1.4, we can decompose $M$ as the 3-sum of matroids $M_1$ and $M_2$ such that $E(M_1)\cap E(M_2)=T$ and $E(M_1)\cap R(M)\neq \phi$. If Lemma 2.3(ii) occurs, $f\in T^*-T'$, then, by Lemma 2.4, $M/f$ has the same regular elements as $M$. By the choice of $M$, $R(M)$ is a circuit-cocircuit of $M/f$. As $R(M)\cup f$ contains a triad of $M$, it follows that $R(M)\cup f$ is not a circuit of $M$. Thus $R(M)$ is a circuit-cocircuit of $M$. 

We may assume by Lemma 2.3 (i) that $M_2$ is regular, $M_1$ is non-regular, and $|R(M)|\subseteq E(M_1)$. By the choice of $M$ if $|E(M_1)|\ge 9$, $R(M)$ is a circuit-cocircuit of ${\rm si}(M_1)$ and therefore of $M$; a contradiction. Thus $M_1$ has at most 8 elements. Since ${\rm si}(M_1)$ is non-regular, ${\rm si}(M_1)$ is isomorphic to $F_7$ or $S_8$. In both cases, $R(M)$ is a circuit-cocircuit of this matroid. $\qed$
\bigskip

Using the previous lemma, we can refine the second part of Lemma 2.4.
\bigskip

\noindent {\bf Lemma 3.2.} {\it Let $M$ be a $3$-connected non-regular binary matroid with $|E(M)|\ge 10$ and suppose $T,T^*$ is an undesired fan of $M$ such that $T^*-T=\{f\}$. Then $M/f$ is a non-regular $3$-connected matroid such that $R(M/f)=R(M)$.}
\bigskip

\noindent {\bf Proof.} We argue by contradiction. Since $T\cap T^*\subseteq R(M)$, it follows from Lemma 2.4 that $|R(M)|\ge 2$.  Lemma 3.1 implies that $|R(M/f)|$ is 2 or 4.  Moreover, $R(M/f)$ is a circuit-cocircuit of $M/f$. So $|R(M)|=2$ and $R(M/f)=4$. 

Since $T^*\subseteq R(M/f)\cup f$, it follows that $R(M/f)$ is also a circuit-cocircuit of $M$. Therefore $T'=T\bigtriangleup R(M/f)$ is a triangle of $M$ and $T'^*=T^*\bigtriangleup R(M/f)$ is a triad of $M$. But $T'$ is a triangle of $M/f$ containing two regular elements of $M/f$ such that no triad of $M/f$ contains these two elements. By Lemma 2.5 these two elements are also regular in $M$. Hence $R(M/f)=R(M)$; a contradiction. $\qed$
\bigskip

A 3-separation $\{X,Y\}$ for a 3-connected matroid is said to be trivial provided $|X|=3$ or $|Y|=3$.
\bigskip

\noindent {\bf Lemma 3.3.} {\it Let $M$ be a $3$-connected non-regular binary matroid such that $|R(M)|\ge 1$. If any non-trivial $3$-separation for $M$ has the union of a triangle and a triad of a undesired fan as one of its sets, then $M$ is isomorphic to $S_8$, $F_7$ or $F_7^*$.}
\bigskip

\noindent {\bf Proof.} If  $|E(M)|\le 8$, then the result holds. Therefore, suppose that $|E(M)|\ge 9$. First assume that $M$ has just one non-trivial 3-separation. By Theorem 1.4, $M$ is the 3-sum of matroids $M_1$ and $M_2$ such that $E(M_1)-E(M_2)$ is the union of the triangle and the triad of the undesired fan. Thus $E(M_1)\cap R(M)\neq \phi$. Observe that Lemma 2.3 (ii) holds in this case. By the uniqueness of the 3-separation for $M$, $M_2$ is internally 4-connected. By Theorem 1.3, $M_2$ is isomorphic to $F_7$. Thus $|E(M)|=8$; a contradiction. Hence $M$ has at least two non-trivial 3-separations.

Let $T_1,T_1^*$ and $T_2,T_2^*$ be different undesired fans of $M$. For $i\in\{1,2\}$, set $Z_i=T_i\cap T_i^*$. By Lemmas 2.4, 3.1, and orthogonality, $R(M)=Z_1\cup Z_2$ is a circuit-cocircuit of $M$. In particular, $Z_1$ and $Z_2$ are unique and these are the unique undesired fans of $M$. If $T_1-T_1^*=\{t\}$ and $T_1^*-T_1=\{f\}$, then $T_2=Z_2\cup t$ and $T_2^*= Z_2\cup f$ because $T_1\bigtriangleup T_2=T_1^*\bigtriangleup T_2^*=R(M)$. Observe that $Z_1\cup Z_2\cup\{f,t\}$ is a 3-separating set for $M$. Thus $|E(M)|=9$ because $M$ has only two non-trivial 3-separations. Hence $M$ is isomorphic to $P_9$ or $P_9^*$; a contradiction. $\qed$

\section{The main result}

In this section we give the proof of Theorem 1.1.
\bigskip

\noindent {\bf Proof of Theorem 1.1.} First, we prove the ``only if'' part. If $M$ is non-binary, then by Theorem 1.2 we may conclude that $M\cong U_{2,4}$. Therefore suppose $M$ is binary and non-regular. Assume that $|R(M)|\ge 2$. If $M$ is an internally 4-connected matroid other than $F_7$ and $F_7^*$, then by Theorem 1.3 $M$ has a minor isomorphic to $M(E_5)$, $S_{10}$, $S_{10}^*$, $T_{12}\backslash e$, or $T_{12}/e$. Observe that $M(E_5)$ and $T_{12}\backslash e$, and $T_{12}/e$ have zero regular elements and $S_{10}$ and $S_{10}^*$ have one regular element; a contradiction because $R(M)\subseteq R(N)$. Thus $M\cong F_7$ and $F_7^*$.
 
We may now assume that $M$ is not internally 4-connected. By Lemma 3.3, $S_8$ is the unique matroid having all non-trivial 3-separations induced by the union of a triangle and a triad of some undesired fan. The result follows in this case. Therefore, we can assume that $M$ has a 3-separation such that none of its sets is the union of a triangle and a triad in a undesired fan, say $\{X_1,X_2\}$. By Theorem 1.4 there are 3-connected matroids (up to parallel elements with the common triangle) $M_1$ and $M_2$ such that $M$ is the 3-sum of $M_1$ and $M_2$ and, for $i\in\{1,2\}$, $E(M_i)=X_i\cup T$. By definition, $T$ is the common triangle between $M_1$ and $M_2$. By Lemma 2.3 we may assume that $M_1$ is non-regular and $M_2$ is regular. Moreover, $R(M)\subseteq X_1$. We may assume that $M_1$ is also 3-connected (the elements in parallel with elements of $T$, if them exist, are in $M_2$) By Lemmas 2.1 and 2.2, $T$ does not span any element of $R(M)$ in $M_1$ or $M_1^*$. Thus by induction we have three possibilities:

First, suppose $M_1$ is isomorphic to $F_7$ or $S_8$. The result follows because $M$ is the 3-sum of a matroid isomorphic to $F_7$ or $S_8$ (that is $M_1$) with a regular matroid (that is $M_2$).

Second, suppose $M_1$ is the 3-sum of matroids $N_1$ and $N_2$ along a triangle $T'$ such that $R(M)\subseteq E(N_1)$; $T'$ does not span any element of $R(M)$ in $N_1$; and $N_1$ is isomorphic to $F_7$ or $S_8$ and $N_2$ is regular (We may assume that $T'\cap E(M_2)=\emptyset$.) If $|E(N_2)\cap T|\ge 2$, then $T\subseteq E(N_2)$ because an element of $E(N_1)-E(N_2)$ spanned by $E(N_2)-E(N_1)$ in $M_1$ must be in parallel with some element of $T'$ in $N_1$. In this subcase, $M$ is the 3-sum of $N_1$ and the regular matroid obtained by doing the 3-sum of $N_2$ and $M_2$ along the triangle $T$. The result follows in this case. Thus we may assume that $|E(N_2)\cap T|\le 1$. As any two triangle of $N_1$ meet (recall that $N_1$ is isomorphic to $F_7$ or $S_8$), it follows that $E(N_2)\cap T=\{t\}$. Thus $t$ is in parallel with ane element $t'$ of $T'$ in $N_2$. Let $N_1'$ be the matroid obtained from $N_1$ by adding $t$ in parallel with $t'$. Note that $T$ is a triangle of $N_1'$. Thus $N_1'$ is isomorphic to $F_7^p$ or $S_8^p$. Moreover, $M$ is the 3-sum of $N_1'$ with $N_2\backslash t$ and $M_2$. The result also follows in this case.

Third, suppose there are matroids $N$, $N_1$, and $N_2$ such that: 
\begin {enumerate}
\item $N$ has elements $t_1$ and $t_2$ in parallel; 
\item $N\backslash t_1$ is isomorphic to $F_7$ or $S_8$; 
\item $E(N_1)$ and $E(N_2)$ are disjoint; 
\item $T_i=E(N)\cap E(N_i)$ is a triangle in both $N$ and $N_i$, for both $i\in\{1,2\}$; 
\item $t_i\in T_i$, for both $i\in\{1,2\}$; 
\item $N_1$ and $N_2$ are regular and 3-connected (up to some parallel elements with elements of $T_1$ and $T_2$ respectively); 
\item $(T_1\cup T_2)\cap E(M_2)=\emptyset$; and 
\item $M_1$ is the 3-sum of $N,N_1$ and $N_2$. 
\end {enumerate}

\noindent We begin by showing that $|E(N_i)\cap T|\le 1$, for both $i\in\{1,2\}$. If $|E(N_i)\cap T|\ge 2$, say $i=2$, then $E(N_2)-T_2$ spans $T$ in $M_1$. As $t_1$ and $t_2$ are the only elements of $N$ in parallel, it follows that $T\subseteq E(N_2)-T_2$, otherwise the unique element belonging to $E(N_2)-T_2$ would be in parallel in $N$ with some element of $T_2$ and this element is not $t_1$. Hence $M$ is the 3-sum of $N,N_1$ and $N_2'$, where $N_2'$ is the 3-sum of $N_2$ and $M_2$ along $T$. The result follows, by induction. Thus we may assume that $|E(N_i)\cap T|\le 1$, for both $i\in\{1,2\}$. Moreover, when $|E(N_i)\cap T|=1$, say $E(N_i)\cap T=\{a_i\}$, $a_i$ is in parallel with some element $a_i'\in T_i$ in $N_i$. If $A_i=\{a_i\}$, when this happens, and $A_i=\emptyset$ otherwise, then $M_1$ is the 3-sum or $N'\backslash[\{a_1,a_2\}-(A_1\cup A_2)]$ with $N_1\backslash A_1$ and $N_2\backslash A_2$, where $N'$ is obtained from $N$ by adding, for both $i\in\{1,2\}$, $a_i$ in parallel with $a_i'$. As $T$ does not span any element of $R(M)$ in $N'$, by Lemma 2.1, and $|R(M)|\ge 2$, it follows that $T$ spans $T_1$ or $T_2$, say $T_2$. That is, each element of $T$ is in parallel with some element of $T_2$ in $N'$. We can transfer these elements for $N_2$ and we arrive at the previous case.

Finally, to see the ``if'' part, we use Lemmas 3.2 and 2.5 to reduce the $S_8$ case to the $F_7$ case in the 3-sums. The $F_7$ case is easy to verify. $\qed$

\bigskip

\noindent {\bf References}

\begin{enumerate}

\item S. R. Kingan, Binary matroids without prisms, prism duals, and cubes, {\it Discrete Mathematics}, {\bf 152} (1996), 211-224.

\item S. R. Kingan, A generalization of a graph result by D. W. Hall, {\it Discrete Mathematics} {\bf 173}, (1997) 129-135. 

\item S. R. Kingan  and M. Lemos, Almost-graphic matroids, {\it Advances in Applied Mathematics}, {\bf 28} (2002), 438 - 477.

\item J.G. Oxley, On nonbinary 3-connected matroids, {\it Trans. Amer. Math. Soc.} {\bf 300} (1987) 663?-679.

\item J. G. Oxley, {\it Matroid Theory}, (1992), Oxford University Press, New York.                                                                                                              
\item P. D. Seymour, Decomposition of regular matroids,  {\it J. Combin. Theory   Ser. B } {\bf 28 } (1980), 305-359.   

\item X. Zhou, On internally 4-connected non-regular binary matroids,  {\it J. Combin. Theory   Ser. B}  {\bf 91 } (204), 327-343.   

\end{enumerate}

\end {document}